\input amstex
\documentstyle{amsppt}
\magnification 1200
\vcorrection{-1cm}

\def\refAL        {1}
\def\refArnold    {2}
\def\refAschen    {3}
\def\refIPS       {4}
\def\refOrevkovFu {5}
\def\refOrevkovHM {6}

\def\Z{\Bbb Z}

\def\Q{\Bbb Q}

\def\LM{\operatorname{LM}}

\def\Hom{\operatorname{Hom}}
\def\End{\operatorname{End}}

\def\condQ   {$(i)$}
\def\condZm  {$(ii)$}

\def\sectMotiv {2}
\def\sectSol   {3}
\def\sectm     {4}
\def\sectQ     {5}

\def\eqbc{1}

\rightheadtext{On modular computation of Gr\"obner bases over $\Z$}

\topmatter

\title     On modular computation of Gr\"obner bases
           with integer coefficients
\endtitle

\author    S.~Yu.~Orevkov
\endauthor
\address   Steklov Math. Institute, Moscow
\endaddress
\address   IMT, Universit\'e Paul Sabatier (Toulouse-3)
\endaddress
\email     orevkov\@math.ups-tlse.fr
\endemail
\endtopmatter

\document

\subhead \S1. Introduction
\endsubhead
Let $A=\Z[X]$, $X=(x_1,\dots,x_n)$. We consider $A$ as a subset of $\Q[X]$,
so, if $I$ is an ideal of $A$, then $\Q I$ is the ideal of $\Q[X]$
generated by $I$. Given a ring $R$, we denote the natural mapping
$A\to A\otimes R=R[X]$ by $\iota_R$.
This note is devoted to the following algorithmic problem
(see [\refAL] for a definition and properties of Gr\"obner bases
of ideals in polynomial rings over $\Z$).

\proclaim{ Problem $(\Cal P)$ }
Suppose that we have an infinite sequence $f_1,f_2,\dots$
of elements of $A$ (a blackbox producing them one by one).
Let $I=(f_1,f_2,\dots)$ be the ideal generated by all the $f_i$'s.
Compute the Gr\"obner base of $I$ under the
assumption that the Gr\"obner base of $\Q I$ and the Gr\"obner bases
of $\iota_{\Z/m\Z}(I)$, $m\in\Z$, are known.
\endproclaim

This problem arises in [\refOrevkovFu], [\refOrevkovHM] (see \S\sectMotiv\ for
more details).
Our solution is based on 
a computation of the exponent of the torsion part of $A/I$ (see [\refAschen]),
and Main Lemma in \S\sectSol\ below.

\subhead\S\sectMotiv. Motivation
\endsubhead
More precisely, the algorithmic problem appeared in
[\refOrevkovFu], [\refOrevkovHM] is the following.
Let $F$ be a finite rank free $A$-module,
$M_0$ a submodule of $F$ (given by a finite generating set $G_0$),
$\tau\in \Hom_A(F,A)$,
and $\rho_1,\dots,\rho_p\in\End_A(F)$.
Compute (a Gr\"obner base of) $\tau(M)$ where $M$ is
the minimal submodule of $F$ such that
$M_0\subset M$ and $\rho_i(M)\subset M$ for any $i=1,\dots,p$.
Theoretically, this problem has an evident solution, namely,
we could compute one by one the Gr\"obner bases of
the modules $M_j$ defined recursively by
$M_{j+1} = \sum_{i=1}^p \rho_i(M_j)$ till the moment when
$M_{j+1}=M_j$. We set at this moment $M=M_j$ and $I=\tau(M)$.

However, if the rank of $F$ is very big, then the computation of Gr\"obner bases of $M_j$'s
may be so long that it cannot be performed in practice
whereas the Gr\"obner bases of $(\Z/m\Z)[X]$-modules $M_j\otimes(\Z/m\Z)$ can be
computed much faster. When we know them, we can compute the Gr\"obner bases of
$\Q[X]$-modules $M_j\otimes\Q$ (see [\refArnold] and \S\sectQ\ below).
Therefore we can compute the modules $M\otimes (\Z/m\Z)$ and $M\otimes\Q$, and hence
the ideals $\iota_{\Z/m\Z}(I)$ and $\Q I$. Also we can successively compute
the generators of $I$ of the form $\tau\rho_{i_1}\rho_{i_2}\dots\rho_{i_k}(g)$, $g\in G_0$.
So, the problem $(\Cal P)$ naturally arises.

\subhead\S\sectSol. Solution to Problem $(\Cal P)$
\endsubhead
Let the notation be as in Introduction.

\proclaim{ Main Lemma }
Let $I$ and $J$ be ideals in $A$ such that $J\subset I$.
Let $m$ be the exponent of the torsion part of the
abelian group $A/J$, i.~e., $m$ is the minimal positive
integer such that $m(\Q J\cap A)\subset J$.
Let $m_1,\dots,m_t$ be pairwise coprime integers whose product is
equal to $m$.
Suppose that
\roster
\item  "\condQ"   $\Q J = \Q I$, i.~e., $I$ and $J$ generate the same ideal in $\Q[X]$,
\item  "\condZm"  $\iota_{\Z/m_i\Z}(J)=\iota_{\Z/m_i\Z}(I)$, $i=1,\dots,t$;
\endroster
Then $I=J$.
\endproclaim

\demo{ Proof } Let $c_i=m/m_i=\prod_{j\ne i}m_j$, $i=1,\dots,t$.
Since $m_i$'s are pairwise coprime, there exist $b_1,\dots,b_t\in\Z$
such that
$$
  b_1c_1+\dots+b_tc_t=1                  \eqno(\eqbc)
$$
We have to show that $I\subset J$.
Let $f\in I$.
By Condition \condZm, for any $i=1,\dots,t$, we have
$f=f_i+m_i h_i$ with $f_i\in J$ and $h_i\in A$.
By Condition \condQ, we have $kf\in J$ for some $k\in\Z$,
hence, for any $i$, we have $km_ih_i = kf-kf_i\in J$, i.~e.,
$h_i$ represents an element
of the torsion part of the group $A/J$.
By definition of $m$, this implies $mh_i\in J$, hence
$c_i f = c_i f_i + c_i m_i h_i = c_i f_i + m h_i \in J$.
Combined with (\eqbc), this yields
$$
   f = b_1(c_1 f) + \dots + b_t(c_tf) \in J
\qed
$$
\enddemo

Thus we obtain the following solution to Problem $(\Cal P)$.
We fix any term order in $A$.
We try one by one 
all ideals $J$ of the form $(f_1,\dots,f_k)$, $k=1,2,\dots$
and for each of them do the following:
\roster
\item
Compute the reduced Gr\"obner bases of the ideals
$\Q I$ and $\Q J$ of the ring $\Q[X]$.
If they do not coincide, then we pass to the next $J$.
\item
Compute $m$ as in Main Lemma using one of the algorithms given
in [\refAschen] (see also \S\sectm).
Let $m=p_1^{\alpha_1}\dots p_t^{\alpha_t}$ be
the factorization of $m$ into prime factors and let $m_i=p_i^{\alpha_i}$.
\item
For each $i=1,\dots,t$,
compute the reduced Gr\"obner bases of $\iota_{\Z/m_i\Z}(I)$ and 
$\iota_{\Z/m_i\Z}(J)$. If they do not coincide, then we
pass to the next $J$. If they do coincide, then we
compute the Gr\"obner base of $J$ and stop the computation.
\endroster

\subhead \S\sectm. Computation of $m$
\endsubhead
Several algorithms are discussed in [\refAschen] to compute
$m$ (defined in Main Lemma).
For the reader's convenience, we reproduce here the algorithm
which (we quote [\refAschen]) `seems better suited for practical
computations'. It consists in the following steps:
\roster
\item
Fix any term ordering and compute a Gr\"obner base $G$ of $J$.
Let $s$ be the least common multiple of the leading coefficients
of the elements of $G$. Then $\Q J\cap A=\Z[1/s]J\cap A$.
\item
Let $Y$ be a new indeterminate. Then
$\Z[1/s]J\cap A=(I,sY-1)\cap A$. A Gr\"obner base $G_s$ of this ideal
can be computed using an elimination term ordering.
\item
Compute $m$ such that $mg\in J$ for any $g\in G_s$.
\endroster

\subhead\S\sectQ. On the modular computation of Gr\"obner bases over $\Bbb Q$
\endsubhead
In our notation, Elizabeth Arnold's theorem can be formulated as follows
(similar theorem holds also for modules).

\proclaim{ Theorem } {\rm[\refArnold; Th.~7.1].}
Let $I$ be a {\bf homogeneous} ideal in $\Z[X]$ and $G\subset\Z[X]$.
Let $\Q I$ be the ideal in $\Q[X]$ generated by $I$.
Let $p$ be a prime number. We set
$I_p=\iota_{\Bbb F_p}(I)$,
$G_p=\iota_{\Bbb F_p}(G)$.
Suppose that:
\roster
\item  $G_p$ is a Gr\"obner base of $I_p$,
\item  $G$ is a Gr\"obner base of the ideal $\langle G\rangle_{\Q}$
       that it generates in $\Q[X]$,
\item  $\Q I\subset\langle G\rangle_{\Q}$,
\item  $\LM(G_p) = \LM(G)$ where ``$\LM$'' denotes the set of leading monomials.
\endroster
Then $G$ is a Gr\"obner base of $\Q I$.
\endproclaim

It is claimed in [\refIPS; Th.~2.4] that the same statement holds for
non-homogeneous ideals. However, it is wrong as one can see in the
simplest example $I=(px+1)$, $G=\{1\}$.

When the ideal is not homogeneous, it can be homogenized by adding a new variable
and after that one can apply E.~Arnold's theorem. For instance, if we homogenize
the above counter-example, i.~e., if we set $I=(px+y)$, $G=\{y\}$, then
the condition (3) no longer holds true, so, the contradiction disappears.

\subhead Acknowledgment
\endsubhead I am grateful to Vladimir Gerdt for useful discussions.

\Refs
\def\r{\ref}

\r\no\refAL
\by W.~W.~Adams, P.~Loustaunau
\book An Introduction to Gr\"obner Bases
\bookinfo Graduate Studies in Mathematics, vol.~3
\publ A.M.S.
\publaddr Providence, RI
\yr 1994
\endref

\r\no\refArnold
\by E.~A.~Arnold
\paper Modular algorithms for computing Gr\"obner bases
\jour J. of Symbolic Computations \vol 35 \yr 2003 \pages 403--419
\endref

\r\no\refAschen
\by M.~Aschenbrenner
\paper Algorithms for computing saturations of ideals
in finitely generated commutative rings.
{\rm Appendix to:}
 Automorphisms mapping a point into a subvariety
\jour J. of Algebraic Geom. \vol 20 \yr 2011 \pages 785--794
(by B.~Poonen)
\endref

\r\no\refIPS
\by V.~Idrees, G.~Pfister, S.~Steidel
\paper Parallelization of modular algorithms
\jour J. of Symbolic Computations \vol 46 \yr 2011 \pages 672--684
\endref

\r\no\refOrevkovFu
\by S.~Yu.~Orevkov
\paper Markov trace on the Funar algebra
\jour arXiv:1206.0765
\endref

\r\no\refOrevkovHM
\by S.~Yu.~Orevkov
\paper Cubic Hecke algebras and invariants of transversal links
\jour Doklady Math. \toappear
\transl arXiv:1307.5862
\endref

\endRefs
\enddocument